\pdfoutput=1
\documentclass[reqno,11pt]{amsart}

\usepackage{amssymb}
\usepackage[margin=1in]{geometry}

\usepackage{tikz}
\usetikzlibrary{calc}

\usepackage{comment}
\usepackage{hyperref}

\newcommand{\vp}{\varphi}

\renewcommand{\b}{\beta}

\renewcommand{\l}{\ell}

\newcommand{\s}{\sigma}
\renewcommand{\t}{\tau}

\newcommand{\w}{\omega}
\newcommand{\x}{\mathbf{x}}
\newcommand{\y}{\mathbf{y}}

\newcommand{\T}{{}^{\mathsf{T}}}
\newcommand{\A}{\mathcal{A}}
\newcommand{\B}{\mathcal{B}}

\newcommand{\F}{\mathbf{F}}

\newcommand{\C}{\mathbb{C}}

\newcommand{\Z}{\mathbb{Z}}
\newcommand{\bs}{\backslash}

\def\NP{\textup{NP}}
\def\P{\textup{P}}
\newcommand{\fig}[1]{Figure~\ref{fig:#1}}
\newcommand{\problem}[1]{\textsc{#1}}
\newcommand{\planar}[1]{\problem{Planar #1}}
\newcommand{\alg}[1]{\textsf{#1}}

\newcommand{\defn}[1]{\textbf{#1}}

\newcommand{\problemdef}[3]{
\medskip
\begin{tabular}{ll}
\multicolumn{2}{l}{\problem{#1}}\\
\textbf{Instance:} & #2 \\
\textbf{Decide:} & #3
\end{tabular}\medskip}

\newtheoremstyle{thm}
  {9pt}{9pt}{\itshape}{}{\bfseries}{}{.5em}{}

\theoremstyle{thm}
\newtheorem{thm}{Theorem}[section]

\newtheoremstyle{defin}
  {9pt}{9pt}{}{}{\bfseries}{}{.5em}{}
\theoremstyle{defin}

\newtheoremstyle{exm}
  {9pt}{9pt}{}{}{\scshape}{}{.5em}{}
\theoremstyle{exm}

\newtheoremstyle{proof}
  {}{}{}{}{\itshape}{:}{.5em}{}
\theoremstyle{proof}

\newcommand{\pakpage}{\newpage}
\renewcommand{\pakpage}{}

\begin{document}
\title{Some NP-complete edge packing and partitioning problems in planar graphs}
\author[Jed~Yang]{Jed~Yang$^\star$}
\keywords{Edge packing, spanning trees, NP-completeness, planar graphs, graph partitioning}
\thanks{\thinspace ${\hspace{-.45ex}}^\star$School of Mathematics, University of Minnesota, Minneapolis, MN 55455, USA; \thinspace \texttt{jedyang@umn.edu}}

\begin{abstract}
Graph packing and partitioning problems have been studied in many contexts,
including from the algorithmic complexity perspective.
%
Consider the packing problem of determining whether a graph contains a spanning tree and a cycle that do not share edges.
Bern\'ath and Kir\'aly proved that this decision problem is \NP-complete
and asked if the same result holds when restricting to planar graphs.
Similarly, they showed that the packing problem with a spanning tree and a path between two distinguished vertices is \NP-complete.
They also established the \NP-completeness of the partitioning problem of determining whether the edge set of a graph can be partitioned into a spanning tree and a (not-necessarily spanning) tree.
We prove that all three problems remain \NP-complete even when restricted to planar graphs.
\end{abstract}
\maketitle

\newcommand{\cf}[1]{\mathsf{#1}}
\renewcommand{\A}{\cf{A}}
\renewcommand{\B}{\cf{B}}
\newcommand{\pf}[1]{\mathsf{#1}}
\newcommand{\Pst}{\pf{P}_{st}}
\newcommand{\SpT}{\pf{SpT}}
\renewcommand{\C}{\pf{C}}
\renewcommand{\T}{\pf{T}}
\newcommand{\Cut}{\pf{Cut}}
\renewcommand{\F}{\pf{F}}
\section{Introduction}
A connected graph contains a spanning tree.
Does it contain two spanning trees which do not share edges?
In other words, can the graph stay connected after removing the edges of a spanning tree?
This problem can be solved in polynomial time \cite{NW,Tutte}.
In general, given two classes $\A$ and $\B$ of graphs,
one could ask the \defn{packing problem} of whether a graph $G$
contains edge-disjoint subgraphs $A\in\A$ and $B\in\B$.
Similarly, one could consider the \defn{covering problem}
where the union of $A$ and $B$ is~$G$,
or the \defn{partitioning problem},
where $A$ and $B$ are edge-disjoint and whose union is~$G$.

Bern\'ath and Kir\'aly~\cite{BK}
considered seven classes,
including paths, cycles, and trees.
They noted that
there are $44$ natural\footnote{For example, it makes no sense to pack paths as a trivial path consisting of a single vertex can be used.} graph theoretic questions
under this setup.
Prior to their work, some of these problems were known to be in~$\P$,
while some were \NP-complete.
They settled the status of each remaining problem in the sense of either giving a polynomial time algorithm
or proving that the problem is \NP-complete.
Moreover, for the \NP-complete problems,
they noted that most of these remain \NP-complete even when restricted to planar graphs.\footnote{Most of these remain \NP-complete even when restricted to planar graphs of maximum degree $3$ or~$4$.}
However, five of these cases were left open.
The goal of this paper is to settle the three remaining cases that involve spanning trees.

\begin{thm} \label{t:main}
The following problems are \NP-complete, even when restricted to planar graphs:
\begin{enumerate}
\item the packing problem with a spanning tree and a cycle;%
\footnote{Previously, this problem was (erroneously) claimed to be in~\P.  See Section~\ref{ss:error} for a discussion.}
\item the packing problem with a spanning tree and a path between two distinguished vertices; and
\item the partitioning problem with a spanning tree and a tree.
\end{enumerate}
\end{thm}

The paper is organized as follows.
In Section~\ref{s:defn}, we give definitions and some background.
The \NP-completeness of the three problems are presented in sections~\ref{s:CSpT}, \ref{s:PstSpT}, and \ref{s:TSpT}, respectively.
We conclude in Section~\ref{s:remark} with some remarks and briefly discuss the remaining two problems (the partitioning problems with two trees and with a cut and a forest).

\pakpage
\section{Definitions and background} \label{s:defn}
Let $G$ be a graph, and write $V(G)$ and $E(G)$ for its vertex and edge sets, respectively.
We are concerned with undirected simple graphs, and follow standard terminology.
Given a (connected) graph $G$, a subgraph $H$ of $G$ is \defn{nonseparating} if $G-E(H)$ is connected.

Following the notation of~\cite{BK}, albeit with a different font,
we write $\C$ for the class of all cycles,
$\SpT$ the class of spanning trees,\footnote{Formally, $\SpT$ depends on (the number of vertices of) the underlying graph in question.  By abuse of notation, we ignore these trivial technicalities.}
$\T$ the class of (not-necessarily spanning) trees,
and $\Pst$ the class of $s$--$t$ paths (paths from $s$ to~$t$), where $s$ and $t$ are distinguished vertices.
The packing problem with classes $\A$ and $\B$ is denoted $\A\wedge\B$,
while the partitioning problem is denoted $\A+\B$.

The decision problem \planar{$\C\wedge\SpT$} takes a \emph{planar} graph $G$ as input,
and outputs whether there is a cycle $Q\in\C$ and a spanning tree $T\in\SpT$ such that $Q$ and $T$ are edge-disjoint subgraphs of~$G$.
As $G-E(Q)$ contains a spanning tree if and only if it is connected,
we are equivalently asking to decide the existence of a \emph{nonseparating} cycle in~$G$.

Theorem~\ref{t:main} asserts that the following three decision problems are \NP-complete.

\problemdef{\planar{$\C\wedge\SpT$}}
{Planar graph~$G$.}
{Does $G$ contain a nonseparating cycle?}

\problemdef{\planar{$\Pst\wedge\SpT$}}
{Planar graph~$G$ and distinguished vertices $s,t\in V(G)$.}
{Does $G$ contain a nonseparating $s$--$t$ path?}

\problemdef{\planar{$\T+\SpT$}}
{Planar graph~$G$.}
{Can $E(G)$ be partitioned into a tree and a spanning tree?}

These three problems are trivially in \NP.
We prove their \NP-hardness by similar reductions from a planar version of boolean satisfiability,
which we define below.

\newcommand{\negate}[1]{\overline{#1}}
\newcommand{\nxi}{\negate{x_i}}
\renewcommand{\nxi}{-x_i}

A \defn{boolean variable} takes \defn{boolean values} $+$ (true) and $-$ (false).
We identify $+$ and $-$ with $1$ and $-1$, respectively,
allowing us to negate and multiply boolean values by inheriting the notations and operations from~$\Z$.
A \defn{literal} is a variable $x$ or its negation~$-x$.
A finite collection of literals is called a \defn{clause}.
A \defn{boolean expression} $\vp$ (in conjunctive normal form)
consists of a set $X=\{x_1,\ldots,x_n\}$ of variables
and a set $C=\{C^1,\ldots,C^m\}$ of clauses.
%
Let the \defn{associated graph} $G_\vp$ be the graph with vertex set $X\sqcup C$ and edges
$$\{x_iC^j:\text{$x_i$ or $\nxi$ occurs in $C^j$}\}\cup\{x_ix_{i+1}:i\in\{1,\dotsc,n\}\},$$
where subscripts are hereafter read modulo~$n$.
A boolean expression $\vp$ is \defn{planar} if its associated graph $G_\vp$ is.
An assignment $f:X\to\{\pm\}$ of boolean values to the variables
is \defn{satisfying} if
each clause contains a $+$ literal under such an assignment.


\problemdef{\planar{SAT}}
{Planar boolean expression $\vp$.}
{Does $\vp$ admit a satisfying assignment?}

Lichtenstein \cite{Lichtenstein} proved that \planar{SAT} is \NP-complete, even when each clause contains precisely three literals.

\pakpage
\section{\planar{$\C\wedge\SpT$} is \NP-complete} \label{s:CSpT}
We reduce \planar{SAT} to \planar{$\C\wedge\SpT$}.
Given a boolean expression $\vp$ with variables $X=\{x_1,\dotsc,x_n\}$ and clauses $\{C^1,\dotsc,C^m\}$
such that the associated graph $G_\vp$ is planar,
we form a new planar graph $H_\vp$ from $G_\vp$ such that 
$H_\vp$ contains a nonseparating cycle if and only if $\vp$ admits a satisfying assignment.

\subsection{Reduction construction} \label{ss:CSpT-construct}
Fix a proper plane drawing of $G_\vp$.
Color each edge $x_iC^j$ by $+$ if $x_i\in C^j$ and $-$ if $\nxi\in C^j$.
(Note that if a clause $C^j$ contains both $x_i$ and $\nxi$,
then we may omit~$C^j$.
Therefore, without loss of generality, assume this does not happen.)

For each $i$, take a small neighborhood $D_i$ about $x_i$ containing only the vertex $x_i$
and initial segments of edges leaving~$x_i$.
We locally modify the plane graph. 
First, subdivide edge $x_ix_{i+1}$ to the path $x_it_is_{i+1}x_{i+1}$,
where the new vertices $s_i$ and $t_i$ lie within~$D_i$ for each~$i$.

From now on, we focus on a fixed~$i$ and perform local replacements.
Each new vertex is added inside~$D_i$,
and is decorated with a subscript~$i$,
which may be suppressed for notational convenience.
Similarly, new edges are to be drawn inside~$D_i$,
the reader is encouraged to check that each new edge can be drawn without introducing crossings.

\begin{figure}[hbtp]
\begin{tikzpicture}
\def \x {1}
\def \y {1}
   \tikzstyle{every node}=[draw,circle,fill=white,minimum size=4pt,inner sep=0pt]

   \draw [ultra thick] (0,0) node (s) [label=left:{$y_{i-1}=s_i$}] {}
      -- ++(\x,+\y) node (u0) [label=above:$u^0$] {}
      -- ++(\x,0) node (u1) [label=above:$u^1$] {}
      -- ++(\x,-\y) node (u2) [label=right:{$u^2=v^2$}] {}
      -- ++(\x,-\y) node (u3) [label=below:$u^3$] {}
      -- ++(\x,0) node (u4) [label=below:$u^4$] {}
      -- ++(\x,0) node (u5) [label=below:$u^5$] {}
      -- ++(\x,+\y) node (u6) [label=right:{$u^6=v^6$}] {}
      -- ++(\x,+\y) node (u7) [label=above:$u^7$] {}
      -- ++(\x,0) node (u8) [label=above:$u^8$] {}
      -- ++(\x,-\y) node (t) [label=right:{$t_i=y_i$}] {}
      ;

   \draw (s)
      -- ++(\x,-\y) node (v0) [label=below:$v^0$] {}
      -- ++(\x,0) node (v1) [label=below:$v^1$] {}
      -- (u2)
      -- ++(\x,+\y) node (v3) [label=above:$v^3$] {}
      -- ++(\x,0) node (v4) [label=above:$v^4$] {}
      -- ++(\x,0) node (v5) [label=above:$v^5$] {}
      -- (u6)
      -- ++(\x,-\y) node (v7) [label=below:$v^7$] {}
      -- ++(\x,0) node (v8) [label=below:$v^8$] {}
      -- (t)
      ;

   \draw (v0) -- ++(0,+\y) node (w0) [label=right:$\w^0$] {} -- (u0);
   \draw (v4) -- ++(0,-\y) node (w4) [label=right:$\w^4$] {} -- (u4);
   \draw (v8) -- ++(0,+\y) node (w8) [label=right:$\w^8$] {} -- (u8);

   \draw (u1) -- ++(\x,0) node (w2) [label=above:$\s^2$] {} -- (v3);
   \draw (v1) -- ++(\x,0) node (z2) [label=below:$\t^2$] {} -- (u3);
   \draw (u5) -- ++(\x,0) node (w6) [label=below:$\s^6$] {} -- (v7);
   \draw (v5) -- ++(\x,0) node (z6) [label=above:$\t^6$] {} -- (u7);
\end{tikzpicture}
\caption{Replacement within $D_i$; $k=2$; $P_i^-$ darkened.}
\label{fig:replacement}
\end{figure}
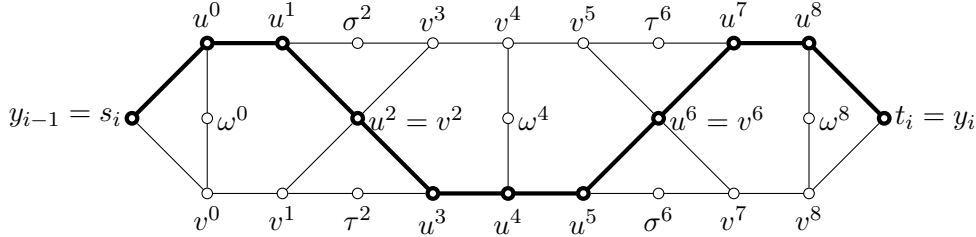

Subdivide each edge $x_iC^j$ with a new vertex~$\b^j=\b^j_i$,
and color the vertex with the color of~$x_iC^j$.
Delete vertex~$x_i$ (and all incident edges).
Let $k$ be a sufficiently large number, say, two times the number of the $\b$ vertices introduced.
Add paths $P_i^+=sv^0v^1\ldots v^{4k}t$ and 
$P_i^-=su^0u^1\ldots u^{4k}t$,
where $u^j=v^j$ for $j\equiv2\pmod4$.
Draw them in such a way that they ``cross'' each other at these common points
(see \fig{replacement}),
and color the edges in $P^+$ and $P^-$ by $+$ and $-$, respectively.
Add path $u^j\w^jv^j$ for each $j\equiv0\pmod4$,
and paths $u^{j-1}\s^jv^{j+1}$ and $v^{j-1}\t^ju^{j+1}$ for each $j\equiv2\pmod4$.

Finally, for each $\b^j$, pick an edge of $P^+$ or $P^-$ of the same color,
subdivide that edge with a vertex~$a^j$,
and join $a^j\b^j$ by an edge.
The edges from the subdivisions stay in $P^+$ or $P^-$ in the obvious way.
As $k$ is large,
we may pick these edges distinct,
such that the new edges $a^j\b^j$ can be properly drawn within $D_i$ without introducing crossings.

After doing such local replacements for all~$i$,
contract $t_is_{i+1}$ to a new vertex $y_i$ for each~$i$.
We consider $P_i^+$ and $P_i^-$ to be $y_{i-1}$--$y_i$ paths.
This concludes the construction of~$H_\vp$.

\pakpage
\subsection{Correctness of reduction} \label{ss:CSpT-proof}
First, note that the construction produces a \emph{planar} graph $H:=H_\vp$ in polynomial time.

Suppose $f:X\to\{+,-\}$ is a satisfying assignment for~$\vp$.
For $r\in\{+,-\}$, let $Q^r=\bigcup_i P_i^{rf(x_i)}$.
(Recall that the multiplication of $\{+,-\}$ is afforded by the identification with $\pm1\in\Z$.)
Note that $Q^+$ and $Q^-$ are edge-disjoint cycles.
We claim that $H-E(Q^-)$ is connected,
which means there is a spanning tree of $H$ that is edge-disjoint from~$Q^-$, as desired.

Indeed, $H-E(Q^-)$ contains the cycle~$Q^+$,
and the $\s$, $\t$, and $\w$ vertices are clearly connected to~$Q^+$.
Moreover, each vertex in $V(Q^-)\bs V(Q^+)$ is connected to $Q^+$ through a $\s$, $\t$, or $\w$ vertex.
As $a^j_i$ and $\b^j_i$ are connected to $C^j$ for each~$j$,
it remains to show that each $C^j$ is connected to~$Q^+$.
By definition, $Q^+$ contains $P_i^{f(x_i)}$.
As $f$ is satisfying, $C^j$ contains a literal $x_i$ such that $f(x_i)=+$, or a literal $\nxi$ such that $f(x_i)=-$.
In either case, $\b^j_i$ is colored with $f(x_i)$,
and hence $a^j_i$ lies on $P_i^{f(x_i)}$.
That is, $C^j\b^j_ia^j_i$ is a path from $C^j$ to $Q^+$, as desired.

\bigskip

Conversely, suppose $H$ contains a cycle $Q$ such that $H-E(Q)$ is connected.
Note that $Q$ cannot contain a vertex of degree $2$ (in $H$),
lest the vertex be isolated in $H-E(Q)$.
Therefore, $Q$ is contained in the subgraph $H'$ of $H$
where the vertices of degree $2$, namely, the $\b$, $\s$, $\t$, and $\w$ vertices, are deleted.
Moreover, the (now) isolated $C$ vertices may also be deleted from~$H'$.
Let $Y$ consists of the $y$ vertices and the $u^j$ vertices for $j\equiv2\pmod4$.
There are two \defn{large} faces that contain $Y$ on its boundary.
Each of the remaining \defn{small} faces is bounded by precisely two vertices in $Y$ and two paths between them.

The $\w^j$ vertices prevent $Q$ from containing (both paths of) a small face,
lest $u^j\w^jv^j$ form a connected component of $H-E(Q)$.
Therefore, $Q$ contains precisely one of the two paths for each small face.
Similarly, the $\s$ and $\t$ vertices force $Q$ to ``oscillate''
and contain either $P_i^+$ or $P_i^-$ between $y_{i-1}$ and $y_i$ for each~$i$.
As such, $Q=\bigcup_i P_i^{-f(x_i)}$ for some $f:X\to\{\pm\}$.

It remains to show that this $f$ is a satisfying assignment.
Suppose not, and there is some clause $C^j$ such that every literal evaluates to $-$ by~$f$.
This means that $a^j_i\in Q$ for each $i$,
and hence $C^j$, together with the $a^j_i$ and the $\b^j_i$ for all~$i$, form a connected component in $H-E(Q)$,
a contradiction.

\pakpage
\section{\planar{$\Pst\wedge\SpT$} is \NP-complete} \label{s:PstSpT}
A similar reduction from \planar{SAT} to \planar{$\Pst\wedge\SpT$} exhibits its \NP-completeness.

\subsection{Reduction construction} \label{ss:PstSpT-construct}
We follow the construction in Section~\ref{ss:CSpT-construct}.
Given a boolean expression $\vp$ with planar~$G_\vp$,
we form $H_\vp$ in the same way,
except that when contracting $t_is_{i+1}$ to a new vertex $y_i$ for each~$i$,
we do not contract $t_ns_1$.
Instead, delete the edge $t_ns_1$, and let $s=s_1$ and $t=t_n$.
Call this graph $H:=H'_\vp$.

\subsection{Correctness of reduction}
Of course, $H$ is a planar graph that is constructed in polynomial time.
From a satisfying assignment,
define $Q^+$ and $Q^-$ the same way
but note that they are edge-disjoint $s$--$t$ paths.
As before, $H-E(Q^-)$ is connected.
Conversely, if there is an $s$--$t$ path $Q$ such that $H-E(Q)$ is connected,
then $Q=\bigcup_iP_i^{-f(x_i)}$, yielding a satisfying assignment~$f$.
We omit the easy details.

\section{\planar{$\T+\SpT$} is \NP-complete} \label{s:TSpT}
We show that \planar{$\T+\SpT$} is \NP-complete by a reduction from \planar{SAT}.

\subsection{Reduction construction} \label{ss:TSpT-construct}
We continue with the construction in Section~\ref{ss:PstSpT-construct}.
Given a boolean expression $\vp$ with planar~$G_\vp$,
we form $H'_\vp$ as above,
but additionally insist that when subdividing edges of $P^+$ and $P^-$ to create the $a$ vertices,
edges must belong to different small faces (as defined in Section~\ref{ss:CSpT-proof}).
Add new cycles $ss's''s$ and $tt't''t$, where $s',s'',t',t''$ are new vertices.
Call this new graph $H:=H''_\vp$.

\subsection{Correctness of reduction}
As before, $H$ is a planar graph constructed in polynomial time.
Suppose there is a tree $T$ such that $H-E(T)$ is a spanning tree.
Certainly $T$ must contain an edge of $ss's''s$ and an edge of $tt't''t$.
As such, $T$ contains an $s$--$t$ path $Q$.
As $H-E(T)$ is a spanning tree, $H-E(Q)$ is connected.
By the same argument as the previous two correctness proofs,
we extract a satisfying assignment~$f$ from~$Q$.

Conversely, take a satisfying assignment $f$ and define $Q^+$ and~$Q^-$ as before.
We alter $Q^+$ and $Q^-$ such that they are a spanning tree and a tree, respectively,
and their edge sets partition that of~$H$.
Take the edges incident to the $\s$, $\t$, and $\w$ vertices,
and add them (along with all incident vertices) to~$Q^+$.
Note that $Q^+$ is still a tree and contains all $u$ and $v$ vertices.

Consider a clause~$C^j$.
Let $i$ be minimal such that $a^j_i\in V(Q^+)$.
As $f$ is satisfying, such an $i$ exists.
Add the path $a^j_i\b^j_iC^j$ to~$Q^+$.
Note that $Q^+$ is still a tree, since we added two new edges and two new vertices.
For each $a^j_\l\in V(Q^+)$, $\l\neq i$,
add $a^j_\l\b^j_\l C^j$ to~$Q^+$.
Since both $a^j_\l$ and $C^j$ were already in~$Q^+$,
we added one new vertex and two edges,
and hence created a single cycle.
The cycle contains an $a^j_\l$--$C^j$ path $P$ avoiding~$\b_\l^j$.
Let $e$ be the first edge in $P$ that intersects $Q^-$,
which exists by the extra stipulation in Section~\ref{ss:TSpT-construct} above.
Delete $e$ from $Q^+$ to destroy the only cycle.
Add $e$ to~$Q^-$, which does not create cycles as they share precisely one vertex.
Finally, for each $a^j_\l\in V(Q^-)$,
add $a^j_\l\b^j_\l C^j$ to~$Q^+$,
which grows the tree $Q^+$ by two vertices and two edges.

After performing the procedure for all clauses,
$Q^+$ is a spanning tree of~$H'_\vp$.
For the remaining $6$ edges of $H-E(H'_\vp)$,
add the paths $ss's''$ and $tt't''$ to $Q^+$ and the remaining two edges $ss''$ and $tt''$ to~$Q^-$.
The edge sets of $Q^+$ and $Q^-$ partition that of~$H$,
while $Q^+$ is a spanning tree and $Q^-$ is a tree,
as desired.

\pakpage
\section{Final remarks} \label{s:remark}
\subsection{}
The use of \planar{SAT} can be replaced by the \NP-complete problem \planar{3SAT},
where each clause has exactly three literals.
If so, the planar graphs constructed here have maximum degree~$4$.
Vertices of degree $4$ are critically used to allow paths to cross each other.
It would be interesting to see if this can be circumvented by using some other ``crossing gadget'' to lower the maximum degree to~$3$.

Similarly, vertices of degree $2$ are used, as in \cite{BK}, to forbid paths or cycles from meandering into the wrong places
and moreover control the way they turn.
If one succeeds in lowering the maximum degree to~$3$,
it would then be reasonable to ask each of these questions when restricted to planar \emph{cubic} graphs, where every vertex is of degree~$3$.

\subsection{} \label{ss:error}
\begin{figure}[hbtp]
\begin{tikzpicture}
\def \x {1}
\def \y {1}
   \tikzstyle{every node}=[draw,circle,fill=white,minimum size=4pt,inner sep=0pt]
   \draw (0,0) node (a) [label=left:$a$] {}
      -- ++(\x,0) node (e) [label=above:$e$] {}
      -- ++(\x,0) node (d) [label=above:$d$] {}
      -- ++(\x,\y) node (b) [label=right:$b$] {}
      ;
   \draw (d) -- ++(\x,-\y) node (c) [label=right:$c$] {};
   \draw (a) -- (b) -- (c) -- (a);
\end{tikzpicture}
\caption{A graph on which the \alg{FindNonSeparatingCycle} algorithm fails.}
\label{fig:error}
\end{figure}
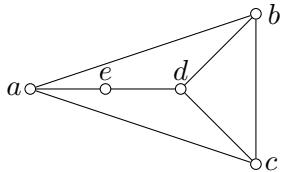

In \cite{BK}, $\C\wedge\SpT$ is listed among problems whose planar restrictions were not known to be \NP-complete.
In an updated arXiv version, \cite{ERROR} is referenced, 
which (erroneously) claims that $\C\wedge\SpT$ is in~$\P$ when restricted to planar graphs.
Indeed, \cite{ERROR} outlines an algorithm \alg{FindNonSeparatingCycle($G$)} that answers the $\C\wedge\SpT$ problem for a planar graph~$G$.
However, the algorithm fails on the graph shown in \fig{error},
which contains a nonseparating cycle $abca$.
If the face bounded by $abdea$ is chosen in Step 1.2 of the algorithm,
the recursive algorithm fails to identify a nonseparating cycle.

\subsection{}
Partition the vertex set of a graph into two non-empty parts.
The set of all edges intersecting both parts form a \defn{cut},
which is \defn{acyclic} if it contains no cycles.
By planar duality, a nonseparating cycle determines an acyclic cut in the dual graph
and, similarly, the existence of an acyclic cut guarantees a nonseparating cycle in the dual graph.
As such, the problem of determining the existence of an acyclic cut in a (planar) graph (possibly with parallel edges) is \NP-complete.

\subsection{}
P\'alv\"olgyi~\cite{Palvolgyi} showed that the problem $\T+\T$ of partitioning a graph into two (not necessarily spanning) trees is \NP-complete.
The reduction is from \problem{NAE-SAT}, where an assignment is satisfying if each clause contains a $+$ and a $-$ (``not all equal'').
The na\"ive approach of simply using the planar version does not work,
since, somewhat surprisingly, \planar{NAE-SAT} can be solved in polynomial time~\cite{Moret}.

The problem $\Cut+\F$ of partitioning a graph into a cut and a forest,
equivalent to coloring the vertex set with two colors such that there are no monochromatic cycles,
is \NP-complete~\cite{BK}.
It is unknown whether these two problems remain \NP-complete when restricted to planar graphs.

\vskip.7cm


\noindent
\textbf{Acknowledgements.} \,
I am grateful to Attila Bern\'ath for the attentive reading of this paper and providing useful comments.
This work is supported by NSF RTG grant NSF/DMS-1148634.

\pakpage
\bibliography{bib}
\bibliographystyle{htam} 

\end{document}